\documentstyle{amsppt}

\NoBlackBoxes

%\textwidth=6.0in
%\textheight=10.0in

%\font\ref=cmr9
%\font\refit=cmti9
%\font\refbf=cmbx9

\nopagenumbers

\input amstex
\topmatter
\title Quantum integrable systems and differential Galois theory
\endtitle

\rightheadtext{Quantum completely integrable systems}

\author
Alexander Braverman, Pavel Etingof
and Dennis Gaitsgory
\endauthor

\address
\newline
A.~B. :  School of Math. Sciences, Tel-Aviv University,
Ramat-Aviv, 69978, Israel.
\newline
D.~G. :  School of Mathematics, IAS, Olden Lane, Princeton, NJ 08540
\newline
P.~E. : Department of Mathematics, Harvard University, Cambridge, MA
02138.
\endaddress

\email
\newline
A.~B. : braval\@ math.tau.ac.il
\newline
P.~E. : etingof\@ math.harvard.edu
\newline
D.~G. : gaitsgde\@ math.ias.edu
\endemail

\abstract
{This paper is devoted to a systematic study of quantum completely
integrable systems 
(i.e. complete systems of commuting 
differential operators) from the point of view of algebraic geometry.
We investigate the eigenvalue problem for such systems and the  
corresponding $D$-module when the eigenvalues are in generic position. 
In particular, we show that the differential Galois group of this 
eigenvalue problem is reductive
at generic eigenvalues. This implies that
a system is algebraically integrable (i.e. 
its eigenvalue problem is explicitly solvable
in quadratures) if and only if the differential Galois group
is commutative for generic eigenvalues. 
We apply this criterion of algebraic integrability to two examples: 
finite-zone potentials and the elliptic Calogero-Moser system.
In the second example, we obtain a proof of the Chalyh-Veselov
conjecture that the Calogero-Moser system 
with integer parameter is algebraically integrable, using
the results of Felder and Varchenko. 
}
\endabstract

\endtopmatter

\def\L{\Lambda}
\def\Spec{\text{Spec}}

\def\Ua{U_q(\tilde\g)}
\def\U2{{\Ua}_2}
\def\g{\frak g}

\def\d{\partial}

\def\l{\lambda}

\def\<{\langle}
\def\>{\rangle}

\define\sbh{\subheading}
\define\Lam{\Lambda}
%\define\wp{{\Cal p}}
\define\OO{{\Cal O}}
\define\Si{S=(\Lam,\theta)}
\define\Sip{S'=(\Lam',\theta')}
\define\MM{{\Cal M}}

\define\Xe{X_{k(\Lam)}}
\redefine\AA{{\Bbb A}}
\define\Ab{{\Cal A}}
\define\Bb{{\Cal B}}
\define\qm{QCIS(X)/\sim}
\define\NN{{\Cal N}}
\define\KK{{\Cal K}}
\define\HH{{\Cal H}}
\define\MMs{{\Cal M}_S}
\define\VVs{{\Cal V}_S}
\define\VV{{\Cal V}}
\define\CC{{\Bbb C}}
\define\kb{Z(\HH)\underset{\HH}\to\otimes k(\HH)}
\define\ot{{\overline t}}

\redefine\d{\partial}
\define\ga{gr(A)}

\heading 0. Introduction 
\endheading
 
\sbh{0.1} Let us recall that in classical mechanics an integrable Hamiltonian 
system on a manifold $X$ of dimension $n$ is a collection of
functions $I_1,...,I_n$ on the cotangent bundle $T^*X$, that Poisson commute 
among themselves and that are functionally independent \cite{Ar}. 

An analogous definition can be given in the framework of algebraic geometry.
Namely, let $X$ be a smooth connected $n$-dimensional algebraic variety
over a field $k$ of zero characteristic which we assume to be algebraically
closed.

By an integrable Hamiltonian system on $X$ we mean
a pair $(\Lambda,f)$, where $\Lambda$ is an irreducible 
$n$-dimensional affine algebraic variety, and $f: T^*X\to\L$ 
is a dominant map whose generic fiber is Lagrangian.
If $\L=\Bbb A^n$, such a mapping is defined by an $n$-tuple
 $I_1,...,I_n$ of regular functions
on $T^*X$ which are algebraically independent and Poisson commute. 
In particular, if $n=1$ ($X$ is a curve), then any nonconstant function 
on $T^*X$ defines an integrable system.  

\sbh{0.2} Now consider the quantum analog of this definition. 
First observe that the map $f$ defines an embedding 
of rings of regular functions
$f^*:\OO(\L)\to \OO(T^*X)$, which is a homomorphism
of algebras with a Poisson commutative image.
The quantum analogue of the Poisson algebra $\OO(T^*X)$ is  
the algebra $D(X)$ of differential operators on $X$.
Thus, by a quantum completely integrable system (QCIS) on $X$ we mean a 
pair $(\Lambda,\theta)$, where $\Lambda$ is an irreducible 
$n$-dimensional affine algebraic variety,
and $\theta: \OO(\L)\to D(X)$ is an embedding of algebras. 
As before, if $\L=\Bbb A^n$, such a mapping is defined by an $n$-tuple 
$D_1,...,D_n$ of differential operators 
on $X$ which are algebraically independent and which commute with one another. 

\sbh{0.3} 
Given a QCIS $S=(\L,\theta)$ on $X$, and a point $\l$ of $\L$,
one may consider a $D$-module $\MM_\lambda$ on $X$ generated by an element
${\bold 1}_\l\in \MM_\lambda$ with the relations 
$$\theta(g)\cdot {\bold 1}_\l=g(\lambda){\bold 1}_\l,\, \forall g\in\OO(\Lam)
.\tag 1$$
Such a $D$-module is the same as a system of differential equations: 
$$
\theta(g)\psi=g(\l)\psi,\, g\in \OO(\L) \tag 2
$$ 
-- the eigenvalue problem
for $S$. 

The rank of $S$ is set to be the dimension of the space
of formal solutions of this system near a generic point of $X$ or
which is the same, the rank of $\MM_\lambda$ at the generic point of $X$.
It is easy to show (Proposition 1.3) that this definition does not
depend on the choice of $\lambda$. We call a QCIS non-degenerate if it
has a non-zero rank. 

The following theorem is one of our main results (cf. 1.5):
\proclaim{Theorem}
Let $\Si$ be a non-degenerate QCIS on $X$ and let $\lambda\in\Lam$
be in a generic position (cf. 1.3). Then
\roster
\item
The $D$-module $\MM_\lambda$ is holonomic and is the Goresky-MacPherson 
extension from an open part of $X$ where it is smooth.
\item
The differential Galois group of $\MM_\lambda$ is a reductive algebraic group. 
\endroster
\endproclaim

\sbh{0.4}
Among QCIS, those having rank $1$ are particularly simple since
the corresponding systems (2) admit an explicit solution (cf. 1.4). The
property of having an explicit solution is shared by a wider class of QCIS's
that are called algebraically integrable (cf. 2.4). By definition, a QCIS 
$\Si$ is said to be algebraically integrable if it is dominated (cf. 2.1, 2.4)
by another QCIS $S'$ with $rk(S')=1$. This is an important notion
that appeared first in the context of soliton theory for $\dim(X)=1$ \cite{Kr}
(cf. 0.5 below).

Theorem 0.3 cited above implies the following criterion for algebraic
integrability of QCIS's in terms of the corresponding differential
Galois groups (cf. 3.4).

\proclaim{Theorem} Let $\Si$ be a QCIS on a variety $X$. The following
conditions are equivalent:
\roster
\item
$S$ is algebraically integrable.
\item
For $\lambda\in\Lam$ in a generic position (cf. 1.3) the differential
Galois group of the equation (2) is commutative. 
\item 
For $\lambda\in\Lam$ in a generic position the differential
Galois group of the equation (2) is triangular (cf. 3.1).
\endroster
\endproclaim

\sbh{0.5}
At the end of the paper, we consider some examples. 
The first example deals with finite-zone ordinary differential
operators \cite{Kr}, which arise in the theory of equations of 
Korteweg-de-Vries type.

Let $X$ be an open subset of 
the affine line ${\Bbb A}^1$ and let $L=d^r/dx^r+...+a_r(x)$
be a differential operator on $X$. Such $L$ is called finite-zone
if there exists another differential operator $L'$, of order prime to $r$,
such that $[L,L']=0$. It is well-known that $L$ is finite-zone if
and only if the QCIS on $X$ defined by $L$ is algebraically integrable.
Therefore, Theorem 0.4 implies the following result (cf. 5.1):

\proclaim{Theorem}
A differential operator $L$ as above is finite zone if and only the 
differential Galois group of the  
equation $L\psi=\l\psi$ is commutative for a generic $\l$. 
\endproclaim

The second example deals with the Calogero-Moser integrable system $CM^n_m$
on an $n$-th power of an elliptic curve $E$ which has a complex
number $m$ as a parameter (cf. 5.3). The system $CM^n_m$ can be obtained 
via the representation theory of affine Kac-Moody algebras \cite{EFK}
and is closely connected to the quantum Hitchin's system in the case
of elliptic curves \cite{BD}. 

We will show that the recent results of Felder and Varchenko who computed
the eigenfunctions of the Calogero-Moser system when $m\in{\Bbb Z}$
imply the algebraic integrability of $CM^n_m$ in this case. This result
was conjectured by Chalyh and Veselov in \cite{CV1}.

{\bf Remark.} 
Although the last result is independent of the non-constructive Theorem 0.4,
the question of writing down explicitly the additional operators that commute
with $\theta(\Lam)$ for $\Si=CM^n_m$ remains open.

\sbh{0.6}
The paper is organized as follows: 

In Section 1 we introduce quantum completely integrable systems and discuss
their basic properties.

In Section 2 we define a partial order on the set of QCIS on a given variety
$X$ and investigate its structure. In 2.4 we study the notion of
algebraic integrability of a QCIS.

In Section 3 we recall basic definitions from the differential Galois theory.
In particular, to a QCIS $\Si$ we attach an algebraic group $DGal(S)$
defined over the field $k(\Lam)$ and study its properties.

In Section 4 we discuss generalizations of the results of the previous
sections to the case when $X$ is a formal polydisc.

In Section 5 we study in some detail the two examples mentioned in 0.5.

In Sections 6 and 7 we present the proofs of the main theorems of the paper.

Some auxiliary algebraic results are proven in the Appendix.

\sbh{0.7 Notations}
Throughout the paper $k$ will be an algebraically closed field of zero
characteristic. If $Y$ is a variety over $k$ and if $K$ is a field
containing $k$ we will denote by 
$Y_K$ the variety $Y\underset {Spec(k)}\to\times \Spec(K)$
over $K$ and by $Y(K)$ the set of $K$-points of $Y$. 
For an irreducible $Y$, we say
that some property holds for almost all points of $Y$, if there exists a
countable collection of proper subvarieties of $Y$ such that 
for any $K\supset k$ this property
holds for all $y\in Y(K)$ that do not belong to these subvarieties.

For an irreducible
affine variety $Y/k$ we denote by $\OO(Y)$ the algebra of regular
functions on $Y$ and by $k(Y)$ the field of rational functions on $Y$. The
Grothendieck generic point of $Y$ is by definition the canonical $k(Y)$-point
corresponding to the embedding $\OO(Y)\subset k(Y)$. 

If $N$ is a quasi-coherent sheaf of $\OO(Y)$-modules on $Y$ and if $y\in Y(K)$,
we denote by $N_y$ the fiber of $N$ at $y$. 

If $Y$ is smooth we denote by $D(Y)$ the ring of differential operators on $Y$.
By definition, a system of (linear) algebraic differential equations on $Y$
is the same as a $D$-module $M$ on $Y$ with a distinguished generator 
$\xi\in M$. (The corresponding system is given
by $D(\psi)=0$ for $D\in D(Y)$ if $D\cdot \xi=0$.) 

For an algebraic group $G$ over $K$ we denote by $Rep(G)$ the category
of its $K$-rational representations and by $Forget$ the forgetful
functor $Rep(G)\to Vect_K$.

\sbh{0.8 Acknowledgments}

We wish to thank J.Bernstein, D.Kazhdan, V.Spiridonov
and A.Varchenko for interesting and stimulating discussions. 
We are also grateful to the anonymous referee whose remarks allowed to 
considerably improve the exposition. The work of P.E. was partially supported
by an NSF postdoctoral fellowship. The work of D.G. was supported by the
NSF grant DMS 9304580.

\heading 1. Quantum completely integrable systems \endheading

\sbh{1.1}
Let $X$ be a smooth connected affine algebraic variety of dimension $n$
over $k$.
\proclaim{Definition}
A quantum Hamiltonian system (QHS) on $X$ is a pair
$S=(\Lam,\theta)$, where $\Lam$ is an irreducible affine algebraic variety
over $k$ and $\theta:\OO(\Lam)\to D(X)$ is an embedding of
algebras.
\endproclaim
The following assertion follows from Bernstein's inequality 
(cf. Corollary A.3(1)):

\proclaim{Lemma}
Let $A\subset D(X)$ be a commutative finitely generated subalgebra.
Then $\dim(A)\leq n$.
\endproclaim

Therefore, if $\Si$ is a QHS on $X$, $\dim(\Lam)\leq n$.
\proclaim{Definition}
A QHS $\Si$ is said to be completely integrable (QCIS) if
$\dim(\Lam)=n$.
\endproclaim

\sbh{1.2}
If $\Si$ is a QCIS, the ring $D(X)$ can be regarded as a left module
over $D(X)\otimes\OO(\Lam)$: $D(X)$ acts on itself by the
left multiplication and $\OO(\Lam)$ acts on $D(X)$ by the right multiplication
via its embedding into $D(X)$.

Let now $K$ be a field containing $k$. If $\lambda\in \Lam(K)$ is a 
$K$-valued point of $\Lam$, we can construct a $D$-module ${\MM}_{\lambda}$
on $X_{K}:=X\underset{Spec(k)}\to\times Spec(K)$ by setting
$${\MM}_{\lambda}:=D(X)\underset{\OO(\Lam)}\to\otimes K.$$

The $D$-module ${\MM}_{\lambda}$ possesses a distinguished generator
${\bold 1}_\l\in D(X)$ which is the image of ${\bold 1}\in D(X)$ 
and it satisfies the relations
$$\theta(g)\cdot {\bold 1}_\l=g(\l){\bold 1}_\l, \forall g\in\OO(\Lam).$$

Therefore, a QCIS $\Si$ can be regarded as a family of $D$-modules on $X$
parameterized by the points of $\Lam$.

In particular, let us consider the Grothendieck generic point (cf. 0.7) of 
$\Lam$. We let $\MMs$ denote the corresponding $D$-module on $\Xe$,
$\MMs:=\simeq D(X)\underset{\OO(\Lam)}\to\otimes k(\Lam)$.
 
\proclaim{Lemma}
\roster
\item
There exists a Zariski dense open subvariety $U_0\subset X\times\Lam$, such 
that $D(X)$ is finitely generated as an 
$\OO(X)\otimes\OO(\Lam)$-module when localized to $U_0$.
\item
For each point $\lambda$ of $\Lam$, the $D$-module ${\MM}_{\lambda}$
is smooth when restricted to $X\times\lambda\cap U_0$.
\item
The $D$-module $\MMs$ on $\Xe$ is holonomic.
\endroster
\endproclaim

The second point of the lemma is an immediate corollary of the first one,
which is proven in A.4. The third point is Corollary A.3(2).

\proclaim{Definition}
The rank $rk(S)$ of $\Si$ is the dimension of the fiber of ${\MM}_S$
at the generic point of $\Xe$.
\endproclaim

It is easy to see that $rk(S)=0$ if and only if the localization of $D(X)$
to the subvariety $U_0$ of Lemma 1.2(1) vanishes.

\sbh{1.3}
Before going further, let us explain how the $D$-module $\MMs$, which 
``lives'' on the variety $\Xe$ over a non-algebraically closed field
$k(\Lam)$ is connected to the $D$-modules ${\MM}_{\lambda}$ obtained
from geometric points $\lambda$ of $\Lam$.

\proclaim{Proposition}
For almost all points $\lambda\in \Lam(K)$ (cf. 0.7) the pair 
$(X_K,{\MM}_{\lambda})$ is a form of $(\Xe,\MMs)$.
\endproclaim

\demo{Proof}
The triple $(X,\Lam,\theta)$ can be defined over a countable algebraically
closed sub-field $k_0$ of $k$. Let ${\OO(\Lam)}_0$ denote the 
ring of functions on $\Lam$ defined over $k_0$. In this case, the 
assertion of the proposition will obviously holds for all those
points in $\Lam(K)$, for which the composition map
${\OO(\Lam)}_0\to\OO(\Lam)\to K$ is an embedding, i.e. for those points
of $\Lambda$ that do not belong to any of the proper subvarieties defined
over $k_0$.
$\square$
\enddemo

\sbh{Remark}
In the sequel points of $\Lam$ that do not belong to any of 
the above subvarieties will be referred to as points in a generic 
position.

\sbh{1.4 First examples}

\sbh{a}
By definition, a QCIS $\Si$ on $X$ with $\Lam\simeq{\AA}^n$ is the same as
a collection of $n$ differential operators $L_i$ on $X$ such
that $L_i$'s pairwise commute and are algebraically independent. 

\sbh{b}
Let $X$ be as above and let $S$ be the pair $(X,i)$, where $i$ is the
natural embedding $\OO(X)\to D(X)$. This is a degenerate QCIS and
its rank equals zero.

\sbh{c}
Let $X\simeq{\AA}^n$ and let $\Si$ be a QCIS on $X$. In this case we can
produce a new QCIS, $Fourier(S)$, by setting $Fourier(S)=(\Lambda,\theta')$,
where $\theta'$ is obtained by composing $\theta$ with the automorphism
of $D({\AA}^n)$ that sends $x_i\to \d/\d x_i$ and 
$\d/\d x_i\to -x_i$. Here $x_i$ are standard coordinates 
on ${\AA}^n$. In particular, the QCIS of the previous example goes
under the Fourier transform to a QCIS generated by the 
$\d/\d x_i$'s. The rank of the latter is $1$ and the
corresponding $D$-module ${\MM}_S$ is identified with $\OO(\Xe)$.

\sbh{1.5 QCIS of rank $1$}
Let us first recall a well-known construction from the theory of $D$-modules.
Let $Y$ be a smooth irreducible algebraic variety. It is easy to see that
there is a bijection between the set of closed rational $1$-forms on $Y$
and the set of isomorphism classes of pairs $(M,\xi)$, where $M$ is
a smooth $D$-module on an open subvariety of $Y$ generated as an $\OO$-module 
by $\xi\in M$ (two such pairs $(M,\xi)$ and $(M',\xi')$ are said to be 
isomorphic if there exists an isomorphism between their restrictions 
to a sufficiently small open subvariety that maps $\xi$ to $\xi'$). 

Explicitly, to a closed $1$-form $\omega$
which is regular on an open subvariety $Y_0\subset Y$ one associates
a $D$-module $M(\omega)$ on $Y_0$ defined as follows: $M(\omega)$ has a
distinguished generator $\xi_\omega\in M(\omega)$ that satisfies the relations:
$$\eta\cdot\xi_\omega=\omega(\eta)\xi_\omega$$ for each vector field $\eta$ on 
$Y_0$. When working over the field of complex numbers this construction
has the following interpretation:

\proclaim{Lemma}
Let $Y_0$ be a variety over ${\Bbb C}$ and let $M$ be a $D$-module
on it generated over $D(Y_0)$ by a section $\xi\in M$. Let also $\omega_i$ be 
closed regular $1$-forms on $Y_0$. The following conditions are equivalent:
\roster
\item
There exists an isomorphism of $D$-modules $M\to \oplus_i M(\omega_i)$ 
that maps $\xi$ to $\oplus \xi_{\omega_i}\in\oplus_i M(\omega_i)$
\item 
The functions $e^{\int {\omega_i}}$ form a basis of the space of local 
solutions of the differential equation defined by the pair $(M,\xi)$ (cf. 0.7).
\endroster
\endproclaim
The proof is straightforward.
\medskip

Let now $\Si$ be a QCIS with $rk(S)=1$ and consider the pair $(\MMs,{\bold 1})$
on $X_{k(\Lam)}$. According to the previous discussion there exists a
rational $1$-form $\omega_S$ on $X_{k(\Lam)}$ that corresponds to
$(\MMs,{\bold 1})$ in the above sense. By definition, there exists an
open subvariety $U_{00}\in X\times \Lam$ such that 
$$\omega_S\in {\Omega}^1(X)\underset{\OO(X)}\to\otimes \OO(U_{00}).$$
We have the following assertion:

\proclaim{Proposition} Let $\Si$ be as above. For each point $\l\in \Lam(K)$ 
the restriction of $\MM_\l$ to $U_{00}\cap X\times\l$ is identified with 
$M(\omega_\l)$, where $\omega_\l$ is the restriction of $\omega_S$
to $U_{00}\cap X\times\l$.
\endproclaim

The proposition follows immediately from the definitions. Moreover, 
we deduce from Lemma 1.5 that when $K={\Bbb C}$, the eigenvalue 
problem 0.3(2) corresponding to $\l$ admits an explicit
solution, namely $e^{\int {\omega_\l}}$.

\sbh{1.6}
Next theorem is one of our main results. Its proof will be given in 
section 7.

\proclaim{Theorem}
The $D$-module $\MMs$ on $\Xe$ is irreducible.
\endproclaim

\proclaim{Corollary}
For almost all points $\lambda$ of $\Lam$, the $D$-module 
${\MM}_{\lambda}$ is semisimple and if $rk(S)\neq 0$ it
is the Goresky-MacPherson extension from its restriction to 
$X\times\lambda\cap U_0$, where it is smooth.
\endproclaim

This corollary is an immediate consequence of the theorem combined with
Proposition 1.3.

\heading 2. Partial order on the set of QCIS's \endheading

\sbh{2.0}
Let $\Si$ and $\Sip$ be two QCIS's on $X$.

\proclaim{2.1 Definition}
\roster
\item
We say that $S'$ dominates $S$ if there exists a map of algebras 
$h:\OO(\Lam)\to\OO(\Lam')$ such that $\theta=\theta'\circ h$.
\item
We say that $S$ and $S'$ are 1-step equivalent if one of them dominates
the other and the corresponding map is a birational equivalence, i.e if
$h$ induces an isomorphism $k(\Lam)\simeq k(\Lam')$.
\endroster
\endproclaim

\proclaim{Proposition}
Let $\Si$ and $\Sip$ be two QCIS on $X$ with 
$S'$ dominating $S$ by means of $h$. Let $\l\in\Lam(K)$ be a $K$-valued point
of $\Lam$ such that the inverse image $h^{-1}(\l)$ of $\l$ consists
of $deg(h)$-many $K$-valued points $\l_i'$ of $\Lam'$. Then
$$\MM_\l\simeq\oplus_{i} \MM_{\l_i'}.$$
\endproclaim
The proof is an immediate corollary of the definitions.

Consider now the set $QCIS(X)$ of all QCIS's on $X$. We introduce an
equivalence relation $\sim$ on $QCIS(X)$ to be the equivalence relation
generated by 1-step equivalences as above. We denote the 
corresponding quotient set by $\qm$.

It follows from Proposition 2.1 that the $D$-modules 
${\MM}_S$ and ${\MM}_{S'}$ on $\Xe$ are canonically isomorphic if $S\sim S'$.
 
The following result is easy to verify.
 
\proclaim{Lemma}
The relation of domination on $QCIS(X)$ descends to a correctly defined partial
order relation $\succ$ on $\qm$. 
\endproclaim

\proclaim{Definition}
\roster
\item
Two elements $S_1,S_2\in\qm$ are called 1-step comparable if
either $S_1\succ S_2$ or $S_2\succ S_1$.
\item
Two elements $S_1$ and $S_2$ are called comparable if they can be connected
by a chain of 1-step comparable elements.
\endroster
\endproclaim  

\sbh{Remark}
It is not difficult to show that two QCIS $S_1$ and $S_2$ are comparable
if and only if there exists a third QCIS $S_3$ with $S_3 \prec S_1,S_2$.

\sbh{2.2}
If $\Si$ is a QCIS on $X$, let $Z(S)$ denote the centralizer of
$\theta(\OO(\Lam))$ in $D(X)$. By definition, this is an associative 
algebra over the commutative algebra $\OO(\Lam)$.

\proclaim{Lemma}
\roster
\item
$Z(S)\underset{\OO(\Lam)}\to\otimes k(\Lam)$ has no zero divisors.
\item
$Z(S)\underset{\OO(\Lam)}\to\otimes k(\Lam)$ is naturally isomorphic to
$End_{D(\Xe)}({\MM}_S)$.
\item
$Z(S)\underset{\OO(\Lam)}\to\otimes k(\Lam)$ is finite-dimensional 
as a $k(\Lam)$-vector space.
\endroster
\endproclaim
\demo{Proof}
The first two points follow readily from the definitions.

To prove the third one, recall (1.2), that ${\MM}_S$ is a holonomic
$D$-module on $\Xe$ and therefore its endomorphism algebra is 
finite-dimensional over the field $k(\Lam)$ \cite{Bo}.
$\square$
\enddemo

The next result is non-obvious and is due to Makar-Limanov [ML]. 
Its proof will be given in section 6 for the convenience of the reader.

\proclaim{Theorem}
The algebra $Z(S)$ is commutative.
\endproclaim

\sbh{Remark}
In the case of $\dim(X)=1$ this result was proven in [Am]. 

\proclaim{Corollary}
$Z(S)\underset{\OO(\Lam)}\to\otimes k(\Lam)$ is
a finite field extension of $k(\Lam)$
\endproclaim

The assertion follows by combining the theorem with Lemma 2.2(1,3). 

\sbh{2.3}
According to 2.1, the algebra $End_{D(\Xe)}({\MM}_S)$ depends only
on the class of $S$ in $\qm$. Therefore, the same holds for
the algebra $Z(S)\underset{\OO(\Lam)}\to\otimes k(\Lam)$.
\proclaim{Proposition}
The set of $S'\in\qm$ with $S'\succ S$ is in a bijection with
the set of $k(\Lam)$-subfields of $Z(S)\underset{\OO(\Lam)}\to\otimes k(\Lam)$.
\endproclaim

\demo{Proof}
If $\Sip$ is a QCIS dominating $\Si$ we assign to it a subfield
in $Z(S)\underset{\OO(\Lam)}\to\otimes k(\Lam)$ which is equal to
$$\OO(\Lam')\underset{\OO(\Lam)}\to\otimes  k(\Lam)\simeq k(\Lam').$$

Vice versa, if $\KK$ is a $k(\Lam)$-subfield
of $Z(S)\underset{\OO(\Lam)}\to\otimes k(\Lam)$, it is easy to see
that one can find a finitely generated algebra $\OO(\Lam')$ over
$\OO(\Lam)$ that maps to $D(X)$ and such that its field of 
fractions maps isomorphically to $\KK$.
$\square$
\enddemo

\proclaim{Corollary 1}
For each element $S\in\qm$ there exists an element $S_{max}\in\qm$ with the
following property:

If $S'\in\qm$ is comparable with $S$, then $S_{max}\succ S'$.
\endproclaim

\demo{Proof}
To prove the corollary, take $S_{max}$ to be the element in $\qm$ 
corresponding to the field $Z(S)\underset{\OO(\Lam)}\to\otimes k(\Lam)$
by the above proposition.
$\square$
\enddemo 

\proclaim{Corollary 2}
Let $S$ be a QCIS with $rk(S)\neq 0$. Then 
$$rk(S)\geq \dim_{k(\Lam)}(Z(S)\underset{\OO(\Lam)}\to\otimes k(\Lam)).$$
\endproclaim

\demo{Proof}
Observe first, that if a QCIS $S$ is dominated by another QCIS $S'$, then
$rk(S)=rk(S')\cdot \dim_{k(\Lam)}(k(\Lam'))$, by Proposition 2.1. Take now 
$S'=S_{max}$. We have: 
$$rk(S)\neq 0\Rightarrow rk(S_{max})\neq 0\Rightarrow rk(S_{max})\geq 1
\Rightarrow rk(S)
\geq \dim_{k(\Lam)}(Z(S)\underset{\OO(\Lam)}\to\otimes k(\Lam)).$$
$\square$
\enddemo

\sbh{Remark}
>From the discussion above one can easily deduce that the algebra $Z(S)$ 
depends only on the comparability class of $S$.

\sbh{2.4 Algebraically integrable systems}

\proclaim{Definition}
A QCIS $S$ is said to be algebraically integrable if it
is dominated by another QCIS $S'$ with $rk(S')=1$.
\endproclaim

\proclaim{Lemma}
\roster
\item
The property of algebraic integrability of a QCIS depends only on
its class in $\qm$.
\item
A QCIS $S$ is algebraically integrable if and only if $rk(S_{max})=1$.
\item
A QCIS $S$ is algebraically integrable if and only if 
$$rk(S)=\dim_{k(\Lam)}(Z(S)\underset{\OO(\Lam)}\to\otimes k(\Lam)).$$
\endroster
\endproclaim

\demo{Proof}
The first point is obvious. 

The second and the third points of the lemma follow as in corollary 2.3(2)
by observing that $rk(S)=rk(S_{max})\cdot 
\dim_{k(\Lam)}(Z(S)\underset{\OO(\Lam)}\to\otimes k(\Lam))$.
$\square$
\enddemo

As before, it is easy to see that algebraic integrability of a QCIS $S$ 
depends only on its comparability class.

\sbh{Remark}
If $\Si$ is an algebraically integrable QCIS, the field 
$Z(S)\underset{\OO(\Lam)}\to\otimes k(\Lam)$ will be called the spectral
field of $S$ and for every $\Sip$ dominating $S$ with $rk(S')=1$,
the variety $\Lam'$ will be called a spectral variety of $S$ (cf. 5.2).

\proclaim{Proposition}
For a QCIS $\Si$ the following conditions are equivalent:
\roster
\item
$S$ is algebraically integrable.
\item
Let $K$ be algebraically closed and let $\l\in\Lam(K)$ be in a generic 
position. 
Then the $D$-module ${\MM}_{\lambda}$ decomposes into a 
direct sum of $rk(S)$ $D$-modules on an open subvariety of $X_K$.
\endroster
\endproclaim

\demo{Proof}
Assume first (2). According to Corollary 1.5, the $D$-module ${\MM}_{\lambda}$
decomposes into a direct sum as above over the whole of $X_K$ and therefore
$\dim_{k(\Lam)}(End(\MM_S))=\dim_K(End({\MM}_{\lambda}))=rk(S)$ and (1)
follows in view of Lemma 2.4(3)

To deduce (2) from (1) consider $S'$ dominating $S$ by means of $h$ with 
$rk(S')=1$ and take $\l\in\Lam(K)$ such that $\Lam'$ is \'etale near $\l$ 
(according to  
Proposition 1.3 the validity of (2) does not depend on a particular choice of 
$\l$ in a generic position). Then (2) follows immediately from Proposition 2.1.
$\square$
\enddemo
  
Assume that for a QCIS $\Si$ the equivalent conditions of Proposition 2.4 
hold and let $\Lam'$ be as in the proof above. It follows from 1.5,
that there exists an open subvariety $U_{00}'\subset X\times\Lam'$
and a $1$-form 
$\omega_{S'}\in {\Omega}^1(X)\underset{\OO(X)}\to\otimes\OO(U_{00}')$
such that for $\l\in\Lam(K)$ as above we have an isomorphism
$${\MM}_{\lambda}\simeq\oplus_i M(\omega_{\l_i})$$
over a sufficiently small open subvariety of $X_K$ (here $\l_i$ are as in
Proposition 2.1). When working over ${\Bbb C}$ we obtain the following
criterion for algebraic integrability of a QCIS:

\proclaim{Criterion-Proposition}
\roster
\item
If $\Si$ is an algebraically integrable QCIS, then for $\l\in\Lam({\Bbb C})$
in a generic position the functions 
$e^{\int {\omega_{\l_i}}}$ for ${\omega_{\l_i}}$ as above form a basis 
for the space of local solutions of the eigenvalue problem 
of 0.3(2) corresponding to the point $\l$.
\item
If for some (or, equivalently, any) $\l\in\Lam({\Bbb C})$ in a generic
position the eigenvalue problem as above
admits a basis of local solutions of the form $ e^{\int {\omega_i}}$,
where $\omega_i$ are some rational $1$-forms on $\Lam_{\Bbb C}$,
the system $S$ is algebraically integrable.
\endroster
\endproclaim

\demo{Proof}
The first assertion follows from the above discussion and from Lemma 1.5.
To prove the second assertion observe, that by Lemma 1.5, the $D$-module
${\MM}_{\lambda}$ is isomorphic to the direct sum $\oplus_i M(\omega_i)$
after a restriction to an open subvariety of $\Lam_{\Bbb C}$ and the
needed result follows from Proposition 2.4.
$\square$
\enddemo

The above criterion is the main tool of proving algebraic integrability of
QCIS's (cf. 5.3).

{\bf Remark.} Existence of a basis of the space of
local solutions of the form 
$ e^{\int {\omega_i}}$ for a special 
(non-generic) $\lambda\in\Lam(\Bbb C)$ does
not imply algebraic integrability. 
For example, one can take any two 1-forms $\omega_1,\omega_2$ on the line 
and write a second order differential operator $L$ 
which annihilates $ e^{\int {\omega_i}}$, $i=1,2$.
The QCIS defined by $L$ has rank 2, and in general it is not algebraically 
integrable. That is, the equation $L\psi=\lambda \psi$ 
does not, in general, have explicit solutions of the form 
$ e^{\int {\omega}}$ for generic $\lambda$, even though it has them
for $\lambda=0$. 

\heading 3. Differential Galois group of a QCIS \endheading

\sbh{3.1}
In this subsection we will recollect certain facts from
the theory of differential Galois groups \cite{Katz},\cite{De}.

Let first $Y$ be a smooth connected algebraic variety over a field $K$ 
of characteristic $0$ (but not necessarily algebraically closed)
and let $y\in Y(K)$ be a $K$-point of $Y$.
Let also $M$ be an $\OO$-coherent $D$-module on $Y$.

We define $C(M)$ to be the minimal full abelian subcategory of the category
of $\OO$-coherent $D$-modules on $Y$ which a) contains $M$,
b) is stable under taking subquotients, 
and c) is stable under tensor products and duality. This is a tensor
category over $K$ endowed with a fiber functor $F_y:C(M)\to Vect_{K}$,
which is defined by taking the fiber at the point $y$:
$\NN\to {\NN}_y$ for $\NN\in C(M)$.

According to \cite{De}, there exists an algebraic group, which
we will denote by $DGal(M,y)$, such that the pair $(Rep(DGal(M,y)),Forget)$
is equivalent to the pair $(C(M),F_y)$. 
The group of $K$-points of $DGal(M,y)$ is identified with
the group of automorphisms of the fiber-functor $F_y$.
The group $DGal(M,y)$ is called the differential Galois group
of $M$ with respect to $y$. 

Since the category $C(M)$ is generated by the object $M$,
the canonical representation of $DGal(M,y)$ on $M_y$ is faithful,
i.e. $DGal(M,y)$ is naturally a subgroup in $GL(M_y)$.

\sbh{Examples}

\sbh{a} By definition, the representation of $DGal(M,y)$ on $M_y$ is
irreducible if and only if $M$ is an irreducible $D$-module on $Y$.

\sbh{b} Let $K$ be algebraically closed. Then $DGal(M,y)$ is a
subgroup of a torus (i.e. is commutative and reductive) if and only
if $M$ decomposes into a direct sum of $\OO$-coherent $D$-modules
of rank $1$.

\sbh{c} If $V$ is a vector space over $K$, we call an algebraic sub-group of 
$GL(V)$ triangular if after passing to an algebraic closure of $K$
it is contained in a Borel subgroup of $GL(V)$.
 
We see that $DGal(M)$ is a triangular subgroup in $GL(M_y)$
if and only if after extending scalars to an algebraic closure of $K$
$M$ admits a filtration whose successive quotients have
rank $1$. When $K={\Bbb C}$, a differential analog of Abel's theorem \cite{Ka}
shows that this condition is equivalent to a possibility to express the 
solutions of the system of differential equations corresponding to $M$ 
(for any choice of a generator $\xi\in M$) through functions and
$1$-forms on $Y$ by taking integrals and exponents of integrals. 
 
\proclaim{Lemma}
Let $K\subset K_1$ be a field extension and let $(Y_1,M_1,y_1)$
be obtained by extension of scalars from $(Y,M,y)$. Then the groups
$DGal(M_1,y_1)$ and $DGal(M,y)\underset{Spec(K)}\to\times Spec(K_1)$
are canonically isomorphic.
\endproclaim 
For the proof the reader is referred to \cite{Katz}.

Let now $M$ be a holonomic D-module on $Y$ (not necessarily $\OO$-coherent).
Assume now that $M$ is $\OO$-coherent in a neighbourhood of $y$.
Then we can define $(C(M),F_y)$ by means of replacing $Y$ by an open 
subvariety $Y_0$ containing $y$ such that $M$ becomes $\OO$-coherent
when restricted to $Y_0$. It is easy to verify, that the pair
$(C(M),F_y)$ constructed in this way does not depend on the choice
of $Y_0$ up to a canonical equivalence. In this case the group $DGal(M,y)$
will be defined as above. In the sequel, $can(M)$ will
denote the canonical object in $C(M)$ corresponding to $M$ 
($M|_{Y_0}$).

Sometimes, when the choice of $y$ 
is fixed we will omit the sub-script $y$ in $DGal(M)$. When $K$
is algebraically closed, differential Galois groups $DGal(M,y_1)$
and $DGal(M,y_2)$ corresponding to different points $y_1$ and
$y_2$ are (non-canonically) isomorphic (cf. \cite{De}). In fact,
one can show using the results of [De] that there exists
a group-scheme ${\Cal G}(M)$ over $Y$ endowed with a flat connection,
whose fiber at a point $y\in Y$ is canonically isomorphic to
$DGal(M,y)$.
  
\sbh{Remark}
Another way to construct the group $DGal(M)$ is via the theory
of differential field extensions, in which situation it is often referred
to as a Picard-Vessiot group \cite{Ka}, \cite{De}.

\sbh{3.2}
Let $\Si$ be a QCIS on $X$ and let 
us choose a point $x_0\in X(k)$ contained in the projection of $U_0$ to $X$ 
(cf. Lemma 1.2). For any $k\subset K$, we will denote again by $x_0$ the 
corresponding $K$-point of $X_K$. For a generic $\lambda\in\Lambda(K)$
the $D$-module ${\MM}_{\lambda}$ on $X_K$ will be $\OO$-coherent in
a Zariski neighborhood of this point due to Lemma 1.2(2).
 
\proclaim{Definition}
The differential Galois group $DGal(S,x_0)$ of a QCIS $S$ is set to be the 
differential Galois group $DGal(\MMs,x_0)$. This is an algebraic group
over the field $k(\Lambda)$.
\endproclaim

As before, when the point $x_0$ is fixed, we will replace
the notation $DGal(S,x_0)$ simply by $DGal(S)$.

By definition, $DGal(S)$ is correctly assigned to $S\in\qm$. By combining
lemma 3.1 with proposition 1.3 we get the following assertion:

\proclaim{Proposition}
For almost all points $\lambda$ of $\Lambda$, the group 
$DGal({\MM}_{\lambda})$ is a form of $DGal(S)$.
\endproclaim 

\sbh{Remark}
The following assertion is unnecessary for us and its proof is left
to the reader:
\proclaim{Lemma}
Let $S'\succ S$ be two elements of $\qm$. There
exists a natural map 
$DGal(S,x_0)\underset{Spec(k(\Lam))}\to\times Spec(k(\Lam'))\to DGal(S',x_0)$
such that
\roster
\item
The map $DGal(S,x_0)\underset{k(\Lam)}\to\times k(\Lam')\to DGal(S',x_0)$
is surjective.
\item
The map $DGal(S,x_0)\to Res^{k(\Lam')}_{k(\Lam)}(DGal(S',x_0))$ is
injective. Here $Res^K_L$ denotes Weil's restriction
of scalars functor.
\endroster
\endproclaim

\sbh{3.3}
Let $S$ be as above. We have the holonomic $D$-module $\MMs$ on $\Xe$.
Let also $can(\MMs)\in C(\MMs)$ be as in 3.1.
Let $\VVs$ be the image of $can(\MMs)$ 
under the equivalence of categories $C(\MMs)\to Rep(DGal(S))$ of 3.1. We have:
$\dim(\VVs)=rk(S)$ and $DGal(S)$ is a closed subgroup of $GL(\VVs)$.

\proclaim{Theorem}
\roster
\item
$\VVs$ is an irreducible representation of $DGal(S)$.
\item
The group $DGal(S)$ is reductive.
\endroster
\endproclaim

\demo{Proof}
By theorem 1.6, $\MMs$ is an irreducible $D$-module on $X_{k(\Lambda)}$,
hence $can(\MMs)$ is an irreducible object of $C(\MMs)$, which implies
the first point of the theorem.

The second point follows from the first one, since the group $DGal(S)$
has been shown to possess a faithful irreducible representation.
$\square$
\enddemo

\proclaim{Corollary}
For almost all points $\lambda$ of $\Lam$ 
the group $DGal({\MM}_{\lambda})$ is reductive.
\endproclaim

It is an interesting question to figure out what reductive groups
can be realized as $DGal({\MM}_{\lambda})$ for some QCIS $\Si$ and for 
a generic $\l\in\Lam(K)$. 

Let us stress that the canonical representation of the group
$DGal({\MM}_{\lambda})$ is not necessarily irreducible.

\proclaim{Proposition}
Let $S$ be an element of $\qm$ with $rk(S)\neq 0$. 
Then $S=S_{max}$ (cf. 2.3) if and only if the 
representation $\VVs$ is absolutely irreducible.
\endproclaim

\demo{Proof}
According to Corollary 2.3(1), $S=S_{max}$ if and only if $End_{D(\Xe)}(\MMs)=
k(\Lam)$, which is in turn equivalent to the condition $End_{DGal(S)}(\VVs)=
k(\Lam)$, by Corollary 1.6. The proof follows now from the fact that a 
representation of a reductive group is absolutely irreducible if and only if 
the algebra of its endomorphisms consists only of scalars.
$\square$
\enddemo 

\sbh{3.4}
Let once again $S$ be a QCIS on $X$.

\proclaim{Theorem} The following conditions
are equivalent:
\roster
\item
$S$ is algebraically integrable. 
\item
$rk(S)\neq 0$ and  the group $DGal(S)$ is a triangular subgroup of $GL(\VVs)$.
\item
$rk(S)\neq 0$ and for a generic point $\l$ of $\Lam$ with values in an 
algebraically closed field, the group $DGal({\MM}_{\lambda})$
is a triangular subgroup of $GL({{\MM}_{\lambda}}_{x_0})$.
\endroster
\endproclaim

\demo{Proof}
Conditions (2) and (3) are equivalent by Proposition 3.1.
In view of Theorem 3.3, condition (3) is equivalent to the fact that 
$DGal({\MM}_{\lambda})$ is a commutative reductive group, which by Example 
3.1(b) is equivalent to the fact that ${\MM}_{\lambda}$ decomposes into a 
direct sum of $rk(S)$ summands. The assertion follows now from Proposition 2.4.
$\square$
\enddemo

\proclaim{Corollary}
$S$ is algebraically integrable if and only if $DGal(S)$ is commutative.
\endproclaim

\demo{Proof}
This follows from the fact that for a reductive group the condition
of being commutative is equivalent to the condition of being a triangular
subgroup of some $GL(n)$.
$\square$
\enddemo

According to Example 3.1(c), when working over ${\Bbb C}$, Theorem 3.4
can be reformulated in the following way: If for some (or equivalently any)
${\Bbb C}$-point $\l$ of $\Lam$ in a generic position the eigenvalue 
problem of 0.3(2) is solvable (in the sense of Example 3.1(c)) then the 
system $S$ is algebraically integrable.

\heading 4. Generalization to the formal case  \endheading

\sbh{4.1}
Here we will generalize results of the previous sections
to the case of integrable systems on a formal polydisc. Therefore in this 
section 
$X$ will denote $Spec(k[[x_1,x_2,...,x_n]])$ and the symbols  
$\OO(X)$, $k(X)$ and $D(X)$ will denote respectively
$k[[x_1,x_2,...,x_n]]$, $k((x_1,x_2,...,x_n))$ and
$\OO(X)[\d/\d x_1,\d/\d x_2,...,\d/\d x_n]$. If $k\subset K$ is
a field extension, we define ${\OO(X)}_K$ (resp. ${k(X)}_K,{D(X)}_K$)
as $\OO(X)\underset k\to\otimes K$ (resp. $k(X)\underset k\to\otimes K,
D(X)\underset k\to\otimes K$). Note that we are taking here algebraic
(non-completed) tensor products.
 
\sbh{4.2} 
\proclaim{Definition}
A quantum completely integrable system (QCIS) on $X$ is a pair
$S=(\Lam,\theta)$, where $\Lam$ is an irreducible affine algebraic variety
of dimension $n$ over $k$ and $\theta:\OO(\Lam)\to D(X)$ is an embedding of
algebras.
\endproclaim

\sbh{Remark}
In the formal setting we will assume that the QCIS's we consider
satisfy the following technical condition: 
$D(X)\underset{\OO(\Lam)}\to\otimes k(\Lam)$ is finitely
generated as an $\OO(X)_{k(\Lam)}$-module. This is done in order
to exclude degenerate QCIS like in Example 1.4(b).

\medskip

As in 1.2, given a QCIS on $X$, 
$D(X)$ can be viewed as a family of $D$-modules parameterized by
$\Lam$. We will consider the specialization of this family at the generic
point $\eta$ of $\Lam$. As in section 1, we will get a module
$\MMs$ over ${D(X)}_{k(\Lambda)}$. By the assumption, $\MMs$ is coherent as 
an ${\OO(X)}_{k(\Lambda)}$-module and $rk(S)$ is defined as before to be 
its rank. 

We have an analog of Theorem 1.6 in the formal setting:
\proclaim{Theorem}
The ${D(X)}_{k(\Lambda)}$-module $\MMs$ is irreducible.
\endproclaim

The proof of this theorem is analogous to the one of Theorem 1.6 
(cf. Section 7). 

\sbh{4.3}
Partial order on the set of QCIS on $X$ is defined as in 2.1
and all the results 2.1-2.3 generalize to the formal case without changes.
In particular, we have:

\proclaim{Theorem}
$Z(S)\underset{\OO(\Lam)}\to\otimes k(\Lam)$ is naturally isomorphic to
$End_{{D(X)}_K}({\MM}_S)$ and is a finite field extension of $k(\Lam)$.
\endproclaim

In fact, the proof of Theorem 2.2 that is given in section 6 applies 
equally to the formal situation.

The definition of algebraic integrability as well as Lemma 2.4 are
restated without any changes. One can formulate analogs of Proposition 2.4
and of Proposition-Criterion 2.4 by means of considering finite extensions
of the field $k(\Lam)$.

\sbh{4.4} 
Let $M$ be an $\OO$-coherent ${D(X)}_K$-module for some $k\subset K$.
As in 2.1 we can define the tensor category $C(M)$ endowed with the
fiber-functor $F_0:C(M)\to Vect_K$ that corresponds to the unique closed
$K$-point of $Spec({\OO(X)}_K)$. The differential Galois group $DGal(M)$
can be defined as in 3.1 and it is an algebraic group over $K$. 

Let us now put $K=k(\Lam)$ and $M=\MMs$ as in 4.2. 
The corresponding differential Galois group will be denoted as before
by $DGal(S)$. This is an algebraic group over $k(\Lam)$ whose category of 
representations is generated by a distinguished object ${\VV}_S$.

\sbh{Remark}
If in the definition of $DGal(M)$ we used the completed tensor product 
instead of ${D(X)}_K$, all differential Galois groups would be trivial.
This happens, in particular, for the $D(X)$-modules ${\MM}_{\lambda}$ 
for $\lambda\in\Lam(k)$. However, the reader can easily check that the group 
$DGal(S)$ is never trivial. 

Theorem 4.2 implies the next result:

\proclaim{Theorem}
The group $DGal(S)$ is reductive. 
\endproclaim

\proclaim{Corollary}
The following conditions are equivalent:
\roster
\item
$S$ is algebraically integrable.
\item
$DGal(S)\subset GL({\VV}_S)$ is triangular.
\item
$DGal(S)$ is commutative.
\endroster
\endproclaim

\heading 5. Examples and applications \endheading

\sbh{5.1}
Let $X$ be a formal line and let $\Lam=Spec(k[\lambda])$. We define
a QCIS on $X$ by setting $\theta(\lambda)=L$, where $L$ is a differential
operator on $X$: 
$$L=a_N(x)\frac{d^N}{dx^N}+...+a_0(X),a_i\in k[[x]], \tag I$$
with $N\geq 1$ and $a_N(0)\neq 0$. We denote this QCIS by $S(L)$,
$rk(S(L))=N$. 

The following definition comes from the soliton theory (cf. \cite{Kr}).

\proclaim{Definition}
The operator $L$ is said to be of algebraic type if there exists
an operator $L'\in D(X)$ such that
\roster
\item
$L'$ commutes with $L$.
\item
$L'$ is regular and semi-simple as an endomorphism of the fiber of 
${\MM}_{S(L)}$ at $0$.
\endroster
\endproclaim
(Here, by a regular semi-simple endomorphism of a vector space 
we mean a linear endomorphism having distinct eigenvalues after passing to an 
algebraically closed field.)

The classification of differential operators of algebraic type is
given in \cite{Kr},\cite{Dr}.

\proclaim{Proposition}
The operator $L$ is of algebraic type if and only if the
QCIS $S(L)$ is algebraically integrable.
\endproclaim

\demo{Proof}
Consider the field $Z(S)\underset{\OO(\Lam)}\to\otimes k(\Lam)$
acting on the fiber of ${\MM}_{S(L)}$ at $0$. The proof of the proposition
follows from the following general fact:

\proclaim{Lemma}
Let $F\subset E$ be a finite field extension and let $V$ be a vector space
over $E$. Then $\dim_E(V)=1$ if and only if there exists $e\in E$
whose action on $V$ is a regular semi-simple automorphism of $V$ when the
latter is viewed as a vector space over $F$. Moreover, in this case $e$
generates $E$.
\endproclaim
$\square$
\enddemo

\sbh{Remark}
Let $L$ be algebraic and let $\Lam'$ be a spectral variety of $S$ (cf. 2.4). 
If $\l_0\in k=\Lam(k)$, the fiber of $\Lam'$ over $\l_0$ bijectively
corresponds to the set of eigenvalues of the operator $L'$ when acting
on the fiber of $\MM_{\l_0}$ at $0$, or which is the same on the space
of the formal solutions of the eigenvalue problem $(L-\l_0)(\psi)=0$. 

By combining Theorem 4.3 with the above proposition we get the following 
assertion:

\proclaim{Theorem}
\roster
\item
If $L$ is of algebraic type, then the group $DGal(S(L))$ is
a commutative reductive group.
\item
If the group $DGal(S(L))\subset GL({\VV}_S)$ is triangular, 
then $L$ is of algebraic type.
\endroster
\endproclaim 

\sbh{Remark}
Algebraic integrability of an operator can be defined in an absolutely
similar way in the geometric situation, i.e. when the formal
line $X$ is replaced by an affine smooth algebraic curve. 

\sbh{5.2 The Lam\'e operator}
As a concrete non-trivial example of an operator of algebraic type,
let us consider the Lam\'e operator. In the next subsection this example
will be generalized to a higher-dimensional case.

Let $k={\Bbb C}$, let $E$ be a smooth elliptic curve, and let 
$\wp$ be the Weierstrass function on $E$. Let also
$\frac{d}{d x}$ be an invariant regular vector field on $E$,
normalized in such a way, that with respect to some local parameter 
$u$ near the origin, $\frac{d}{d x}=\frac{d}{d u}+O(u)$,
$\wp(x,E)=u^{-2}+O(u^{-1})$.    
Let $m\in k$, and let
$$
L_m=\frac{d^2}{dx^2}-m(m+1)\wp.\tag II
$$
This operator is called the Lam\'e operator \cite{WW}.
It defines a quantum completely integrable system on $X=E\setminus 0$ with $\Lam={\AA}^1$
as in 5.1. Eigenfunctions of the Lam\'e operator are called Lam\'e functions.

Introduce an anti-involution $*: D(X)\to D(X)$
by $\frac{d}{dx}^*=-\frac{d}{dx}$, $f^*=f$, $f\in\OO(X)$. 
Then the operator $L_m$ is formally selfadjoint, i.e.  
$L_m^*=L_m$.

\proclaim{Proposition} (see \cite{CV1,CV2} and references therein)

The operator $L_m$ is of algebraic type (=finite zone) if and only if 
$m$ is a non-negative integer.
For such $m$ there exists a unique operator $Q_m$ of order $2m+1$
and highest coefficient $1$ such that $[Q_m,L_m]=0$, and $Q_m^*=-Q_m$.
It satisfies an equation $Q_m^2=P_m(L_m)$, where
$P_m\in k[t]$ is a monic polynomial of degree $2m+1$.
\endproclaim 

Let us now show how to reprove the fact that the 
Lam\'e operator is of algebraic type using the results of 2.4. For this purpose
let us recall a construction of Lam\'e functions due to Hermite.

\proclaim{Lemma}
\roster
\item 
Let $f$ be a rational function on $E$. Then the following conditions
are equivalent:

(a) $\pi(f):=f^2+\frac{df}{dx}-m(m+1)\wp$ is a constant function.

(b) $f$ has poles only at $m+1$ points $a_0,a_1,..,a_n$ with
$a_i\neq a_j$ for $i\neq j$ and $a_0=0$ such that 
$(Res_{a_i}(fdx)=1$ for $i\neq 0$,
$(Res_{a_o}(fdx)=-m$ and
$f(a_i+x)+f(a_i-x)$ has a zero at $x=0$ for $i\neq 0$.
\item
The set $H_m$ of functions $f$ satisfying (a) or (b) above has 
a natural structure of a hyperelliptic curve of genus $m$ mapping
to ${\AA}^1$ by means of $f\to\pi(f)$ with the involution $\sigma$ defined by:
$\sigma(f)(x)=f(-x)$.
\endroster
\endproclaim

If now we put 
$\psi_f(x)=e^{\int f(x)dx}$ for $f\in H_m$
we will get a Lam\'e function corresponding to the eigenvalue $\pi(f)$. 
Therefore for any $\l\in {\AA}^1$ which is a non-critical value 
of the projection $\pi$ we get a basis to the space of solutions of
the eigenvalue problem $(L_m-\l)(\psi)=0$ of the form $e^{\int \omega}$
and the algebraic integrability of $S(L_m)$ readily follows from 
Proposition-Criterion 2.4(2). 

\sbh{Remark}
Note that $H_m$ is identified with a spectral variety of $S(L_m)$.

\sbh{5.3 Calogero-Moser systems  and the Chalyh-Veselov conjecture}

Let $E$ be an elliptic curve as above. 
The Calogero-Moser system is a quantum completely integrable
system on the variety $X=E^n\setminus D$, where $D$ is the divisor of 
diagonals.

Let $\frac{d}{dx_i}$ be the vector fields on $X$ defined as in 5.2. 
The following result is due to Olshanetsky and Perelomov \cite{OP}:

\proclaim{Proposition} There exist pairwise
commuting $S_n$-invariant differential operators $L^1...L^n$ in $D(X)$
such that 
$$
\align
&L^1=\sum_{i=1}^n\frac{d}{dx_i},\\
&L^2=\sum_{i=1}^n\frac{d^2}{dx_i^2}-m(m+1)\sum_{i\ne j}\wp(x_i-x_j),\\ 
&L^j=\sum_{i=1}^n\frac{d^j}{dx_i^j}+\text{ lower order terms }.
\tag III
\endalign
$$
The operator $L^j$ is defined by these conditions uniquely up to 
adding a polynomial of $L^1,...,L^{j-1}$. 
\endproclaim

The operators $L^1,...,L^n$ define a quantum completely integrable system on $X$ which
we denote by $CM^n_m$. We have $CM^n_m=(\Lam,\theta)$, where 
$\Lam ={\AA}^n=Spec(\CC[t_1,..,t_n])$, $\theta(t_i)=L_i$.  

The rank of $CM^n_m$ equals $n!$ and for
$n=2$ the system (III) reduces to the Lam\'e equation, which
we considered in the previous section.

\proclaim{Theorem} (a conjecture from \cite{CV1})
Let $m$ be a non-negative integer.
Then the system $CM^n_m$ is algebraically integrable.
\endproclaim

\demo{Proof}
To prove this theorem it is enough to show that for almost all 
$\lambda\in{\AA}^n$, the eigenvalue problem as in 0.3(2) for $CM^n_m$
admits a basis of local solutions of the form $e^{\int \omega}$, where
$\omega$ are some rational closed $1$-forms on $X$.

Solutions of this system of differential equations were recently obtained
by \cite{FV}:

In \cite{FV}, Section 7, the authors construct an irreducible
affine algebraic variety $H_m^n$, which is described as follows. 

Let $e_i$ be the standard basis of $k^n$, ${\frak h}$ be the hyperplane
in $k^n$ defined by the equation $\sum_i y_i=0$, $\alpha_i=e_i-e_{i+1}$, 
$i=1,...,n-1$. Let
$R_m^n$ be the set of rational 
functions $f: E\to {\frak h}$ such that $f$ has a simple pole with residue
$-m((n-1)\alpha_1+...+\alpha_{n-1})$ at the origin and 
$mn(n-1)/2$ other simple poles
$a_1,...,a_{m(n-1)}$:
$m(n-1)$ poles with residue $\alpha_1$, $m(n-2)$ poles
with residue $\alpha_2$,...,$m$ poles with residue $\alpha_{n-1}$. 
Let $H_m^n$ be the set of those $f\in R_m^n$ for which the scalar-valued 
function $g_i(x)=<f(a_i+x)+f(a_i-x),\text{Res}_{a_i}f>$ 
has a zero at $x=0$ for all $i$ (here $<,>$ denotes the usual inner product 
in $k^n$). 
These equations are called the Bethe ansatz equations, and the set
$H_m^n$ has a natural structure of an algebraic variety
called the Hermite-Bethe variety. The dimension of $H_m^n$ 
equals $n-1$. In the case $n=2$, this situation reduces to the one described
in the previous section.
 
The following assertion follows from theorem 7.1 of \cite{FV}:
\proclaim{Lemma} 
There exists a regular map $\pi: H_m^n\times\Bbb A^1\to\Bbb A^n$, 
such that the generic fiber of $\pi$ contains $n!$ points, and a family 
$\Phi(x,f)$ of closed rational 1-forms on $X$ parameterized by $f\in H_m^n$, 
such that 

(i) the formal function
$$
\psi_{f,t}(x)=e^{t\langle x,e_1+...+e_n\rangle +\int \Phi(x,f)},
$$
regarded as a function of $x$, satisfies the differential equations
$$
L_i\psi_{f,t}=\pi(f,t)_i\psi_{f,t},
$$
and

(ii) the functions $\psi_{f,t}(x)$, where $(f,t)$ runs over 
$\pi^{-1}(\lambda)$, are 
linearly independent for a generic $\l$.
\endproclaim

The assertion of algebraic integrability of $CM^n_m$ follows now
from Proposition-Criterion 2.4. 
$\square$
\enddemo

\sbh{Remark}
The proof presented above is rather inexplicit: from the algebraic 
integrability of $CM^n_m$ it follows that $H^n_m\times {\AA}^1$ is 
a spectral variety of $CM^n_m$. This means that there
exists a big sub-algebra inside $\OO(H^n_m\times {\AA}^1)$ mapping
to $D(X)$ and such that its image commutes with $\theta(\OO(\Lam))$.
However, we do not know how to write it down explicitly.

\heading 6. Proof of the commutativity theorem \endheading

\sbh{6.0}
Our objective here is to prove theorem 2.2. 
In this section, we will consider objects endowed with an increasing
filtration. If $C$ is a filtered algebra, $C_n$ will denote the
$n$-th term of its filtration and if $t\in C_n$, ${\sigma}_n(t)$
will denote the image of $t$ in $C_n/C_{n-1}$.

By a filtered map between two objects we will mean
a map strictly compatible with filtrations. In particular, if $C_1\to C_2$
is a filtered embedding of algebras, then $gr(C_1)\to gr(C_2)$ is
an embedding as well.

\sbh{6.1}
Let $X$ be a smooth connected affine variety 
and let $\HH\subset D(X)$ be a commutative subalgebra.
Let $Z(\HH)$ denote the centralizer of $\HH$ in $D(X)$. 
Assume that $\kb$ is finite-dimensional as a vector space
over $k(\HH)$. We want to prove the following more general assertion:

\proclaim{Theorem}
$Z(\HH)$ is commutative as well.
\endproclaim

Theorem 2.2 follows from Theorem 6.1 in view of Lemma 2.2(3).

\sbh{6.2}
Let $R'$ be an associative algebra without zero-divisors
and let $R\subset R'$ be a central subalgebra.

\proclaim{Definition}
\roster
\item
We say that an element $r'\in R'$ is algebraic over $R$ if
there exists a non-zero polynomial $p(x)\in R[x]$ such that $p(r')=0$. 
\item
We say that $R'$ is algebraic over $R$ if every element $r'\in R'$
is.
\endroster
\endproclaim

\proclaim{Proposition}
Let $\Ab\subset \Bb$ be an embedding of filtered associative algebras  
such that the following conditions hold:
\roster
\item
$\Ab$ is central in $\Bb$
\item
$\Bb$ is algebraic over $\Ab$
\item
$gr(\Bb)$ is commutative and has no zero-divisors.
\endroster
Then $\Bb$ is commutative as well.
\endproclaim

Let us first deduce Theorem 6.1 from the above proposition.

\demo{Proof of Theorem 6.1}
We must check that the three conditions above hold when we put
$\Ab=\HH$, $\Bb=Z(\HH)$.

We define the filtrations on $\Ab$ and on $\Bb$ as being induced from
the natural filtration on $D(X)$. By definition, $\Ab\subset \Bb$ is
a filtered embedding and (1) and (3) hold also by definition. 

(2) holds since $\kb$ is finite over $k(\HH)$.

Therefore, $Z(\HH)$ is commutative, by Proposition 6.2.
$\square$
\enddemo

\sbh{6.3}

Proposition 6.2 will follow from a series of auxiliary
statements.

\proclaim{Lemma 1}
Let $\Ab'\subset \Bb'$ be an embedding of commutative Poisson algebras,
such that the following conditions hold:
\roster
\item
$\Bb'$ (and hence $\Ab'$) has no zero-divisors.
\item
$\Bb'$ is algebraic over $\Ab'$
\item
$\Ab'$ Poisson-commutes with $\Bb'$
\endroster
Then the Poisson bracket on $\Bb'$ is zero.
\endproclaim

\demo{Proof}
Let $\ot$ and $\ot'$ be two elements in $\Bb'$. 
We have to show that $\{\ot,\ot'\}=0$.
Take $p\in \Ab'[x]$ to be a non-zero polynomial of minimal degree such that
$p(\ot)=0$. We have:
$$o=\{p(\ot),\ot'\}=\{\ot,\ot'\}\cdot p'(\ot).$$

Since the ring $\Bb'$ has no zero-divisors, $\{\ot,\ot'\}\neq 0$ would imply 
that $p'(\ot)=0$, which is a contradiction, since $deg(p')=deg(p)-1$.
$\square$
\enddemo{QED} 

\proclaim{Lemma 2}
Let $\Ab\subset \Bb$ be an embedding of filtered 
satisfying the conditions of 6.2.
Then $\Bb':=gr(\Bb)$ is algebraic over $\Ab':=gr(\Ab)$.
\endproclaim

\demo{Proof}

Since $\Bb'$ is commutative the set of elements of $\Bb'$ algebraic
over $\Ab'$ is a subalgebra. Therefore, to prove the lemma
it suffices to check that every homogeneous element in $\Bb'$
satisfies a polynomial equation over $\Ab'$.

Let $\ot$ be a homogeneous element of degree $n$ in $\Bb'$. Let
$t$ be an element in $\Bb$ with ${\sigma}_n(t)=\ot$.

Choose a polynomial 
$$p[x]=a_m\cdot x^{m}+a_{m-1}\cdot x^{m-1}+\dots+a_1\cdot x+a_0,$$
$a_i\in\Ab$  with $p(t)=0$. Put $k=max(deg(a_i)+ni)$ and consider
$$q[x]:=\underset{deg(a_j)+nj=k}\to\Sigma 
{\sigma}_{deg(a_j)}(a_j)\cdot x^{j}$$.

We have ${\sigma}_k(p(t))=q(\ot)=0$, and this proves
our assertion.
$\square$
\enddemo

\demo{Proof of Proposition 6.2}
Suppose on the contrary that $\Bb$ were not commutative. Let $k$
be the maximal integer, for which $t\cdot t'-t'\cdot t\in {\Bb}_{i+j-k}$
for every $i,j$ and $t\in {\Bb}_i,t'\in {\Bb}_j$.

Then $\Bb':=gr(\Bb)$ acquires a Poisson structure in a standard manner:
for $\ot$ and $\ot'$ homogeneous elements of degrees $i$ and $j$ respectively,
set $[\ot,\ot']=\sigma_{i+j-k}(t\cdot t'-t'\cdot t)$, where $t$ and $t'$
are any two elements in ${\Bb}_i$ (resp. ${\Bb}_j$) with
$\sigma_{i}(t)=\ot$ (resp. $\sigma_{j}(t)=\ot$).

This is a Poisson structure on $\Bb'$ that satisfies the conditions
of lemma 1. Hence it is trivial and this contradicts 
the assumption of maximality of $k$.
$\square$
\enddemo

\heading 7. The irreducibility theorem \endheading

\sbh{7.1}
Let $X$ be a smooth connected algebraic variety over a field $k$ of
dimension $n$ and let $S=(\Lam,h)$ be a QCIS on $X$. 
Here we will present a proof of Theorem 1.6.

\proclaim{Theorem 1.6}
The $D$-module $\MMs$ on $\Xe$ is irreducible.
\endproclaim

\sbh{7.2}
\demo{Proof}
Let us suppose on the contrary that $\MMs$
is reducible. In such a case there exists a short exact sequence of 
$D$-modules on $\Xe$:
$$0\to {\NN}_1\to\MMs\to {\NN}_2\to 0$$
with $\NN_i\neq 0$ for $i=1,2$.

Put now $N_1:={\NN}_1\cap D(X)\subset\MMs\simeq D(X)\underset{\OO(\Lam)}
\to\otimes k(\Lam)$ and $N_2:=D(X)/N_1$.

The following lemma is straightforward:
\proclaim{Lemma 1}
$N_1$ is a $D(X)\otimes\OO(\Lam)$-submodule of $D(X)$ and
$N_i\underset{\OO(\Lam)}\to\otimes k(\Lam)\simeq{\NN}_i$
for $i=1,2$.
\endproclaim

Consider now $N_2$ as a $D$-module on $X$. 
\proclaim{Lemma 2}
$\dim_{D(X)}(N_2)=2n$
\endproclaim
\demo{Proof of lemma 2}
On the one hand, we have:
$$
\align
&2n\geq \dim_{D(X)}(N_2)\overset{A.3}\to \geq \dim_{D(X)\otimes\OO(\Lam)}(N_2)
\overset{A.2}\to\geq \\
&\dim_{D(\Xe)}({\NN}_2)+\dim(\Lam)=\dim_{D(\Xe)}({\NN}_2)+n
\endalign
$$
On the other hand, by Bernstein's inequality, we have 
$\dim_{D(\Xe)}({\NN}_2)\geq n$
and therefore $\dim_{D(X)}(N_2)=2n$.
$\square$
\enddemo

Now we are ready to conclude the proof of the theorem
by making the following observation:

\proclaim{Lemma 3}
Let $X$ be a smooth connected algebraic
variety and let $M$ be a non-zero $D$-module on $X$
which is a quotient of $D(X)$. Assume that $\dim_{D(X)}(M)=2\dim(X)$.
Then $M=D(X)$.
\endproclaim

By applying this lemma to $N_2$ we see that
either $N_1$ or $N_2$ is zero, which is a contradiction. Therefore the
theorem is proved.

$\square$
\enddemo

\heading Appendix  
\endheading

\sbh{A.0}
We will prove here several statements mainly from commutative algebra
whose proof was omitted in the main body of the paper. We will state
our results in the geometric case (when $X$ is a usual variety). 
The case of a formal polydisc is treated similarly. 

\sbh{A.1}
Let $A$ be an associative filtered algebra $A=\cup A_i$, $i\geq 0$, 
such that the associated graded
algebra $\ga$ is a finitely generated commutative algebra over a field $k$. 
If $M$ is a finitely generated $A$-module, we define $\dim_{A}(M)$ in the 
following way: choose a good filtration
on $M$ (i.e. a filtration compatible with the filtration on $A$ and
such that $gr(M)$ is finitely generated over $\ga$) and set $\dim_{A}(M)$
to be the dimension of the support of $gr(M)$ inside $Spec(gr(A))$.
When $M=0$, $\dim_{A}(M):=-\infty$.

The following fact is well-known:
\proclaim{Lemma}
Let $M$ be an $A$-module and let $\alpha: M\to M$ be an injective endomorphism.
Then $\dim_{A}(M)\geq \dim_{A}(M/Im(\alpha))+1$.
\endproclaim

\sbh{A.2}
Let $A$ be as above and let $A'$ be of the form
$A\otimes C$, where $C$ is a commutative finitely generated domain.
We endow $A'$ with a filtration: $A'_i=A_i\otimes C$. Let $M$
be a finitely generated $A'$-module and let $k(C)$ denote the field of
fractions of $C$.

\proclaim{Lemma}
$\dim_{A'}(M)\geq \dim_{A\otimes k(C)}(M\underset C\to\otimes k(C))+\dim(C)$
\endproclaim

The proof is easy.

\sbh{A.3}
Let $A$, $A'$ and $M$ be as above and let us assume moreover that
$M$ is finitely generated as an $A$-module.

\proclaim{Proposition}
$\dim_A(M)\geq \dim_{A'}(M)$
\endproclaim

\demo{Proof}
Without restriction of generality, we may assume that $C=k[T]$, i.e. $A'=A[T]$.

Consider the $A[T]$-module $M[T]:=M\underset{A}\to\otimes A[T]$.
It is easy to see that $\dim_{A[T]}(M[T])=\dim_A(M)+1$.

We have a short exact sequence of $A[T]$-modules:
$$0\to M[T]\overset{\alpha}\to\to M[T]\overset{\beta}\to\to M\to 0,$$
where $M$ is endowed with the initial $A\simeq A[T]$-structure, 
and the maps $\alpha$ and $\beta$ are defined as follows: 

$$\alpha(m\otimes 1)=T\cdot m\otimes 1-m\otimes T,\
\beta(m\otimes 1)=m.$$

The proposition now follows from Lemma A.1.
$\square$
\enddemo

\proclaim{Corollary}
\roster
\item
Let $S$ be a QHS on a variety $X$. Then $n\geq \dim(\Lam)$.
\item 
If $S$ is a QCIS, then $\MMs$ is holonomic.
\endroster 
\endproclaim

\demo{Proof}
To prove the first point let us apply proposition A.3 to the case when 
$A=D(X)$, $C=\OO(\Lam)$ and $M=D(X)$. We have: 
$$
\align
2n\geq \dim_{D(X)}(D(X))\overset{A.3}\to\geq \dim_{D(X)\otimes\OO(\Lam)}(D(X))
\overset{A.2}\to\geq \\
&\dim_{D(\Xe)}(\MMs)+\dim(\Lam)
\endalign
$$
However, by Bernstein's inequality we have $\dim_{D(\Xe)}(\MMs)\geq n$
($\MMs\neq 0$) and we get: $2n\geq n+\dim(\Lam)$ hence the assertion.

If now $\dim(\Lam)=n$, we necessarily have $\dim_{D(\Xe)}(\MMs)=n$
$\square$ 
\enddemo

\sbh{A.4}
\proclaim{Lemma 1}
Let $A$ and $B'$ be two commutative finitely generated algebras
with $A$ being an integral domain. Let $\phi:A\to B'$ 
be an algebra map such that $B'$ finitely generated over $A$
as an algebra. Let also $M$ be a $B'$-module. Then
\roster
\item
There exists a Zariski dense open subset $U'\subset Spec(A)$ such that
$M$ is flat when localized to $U'$.
\item
Assume moreover that $\dim_{B'}(M)=\dim(A)$. Then there exists a Zariski dense
open subset $U''\subset Spec(A)$ such that $M$ is finitely generated 
over $A$ when localized to $U''$.
\endroster
\endproclaim

Next lemma is an immediate corollary of the one above.

\proclaim{Lemma 2}
Let $A$ be as above and let $B=\cup_{i\in{\Bbb Z}_+}B_i$ be a filtered 
associative algebra with $A$ mapping to $B_0$. Assume moreover, that
$B':=gr(B)$ is commutative and finitely generated as an algebra over $A$.
Let $M$ be a finitely generated $B$-module. Then:
\roster
\item
There exists a Zariski dense open subset $U'\subset Spec(A)$ such that
$M$ is flat when localized to $U'$.
\item
Assume moreover that $\dim_B(M)=\dim(A)$. Then there exists a Zariski dense open
subset $U''\subset Spec(A)$ such that $M$ is finitely generated 
over $A$ when localized to $U''$.
\endroster
\endproclaim

Lemma 1.2(1) follows now from Lemma A.4(2) applied to $A=\OO(X\times\Lam)$,
$B=D(X)\otimes \OO(\Lam)$ and $M=D(X)$.

\enddocument